\documentclass[11pt,twoside]{article}
\usepackage[left=3cm,right=3cm,top=3cm,bottom=3cm]{geometry}
\usepackage{amssymb}
\usepackage{graphicx}
\usepackage{graphicx,type1cm,eso-pic,color}
\usepackage{amssymb}
\usepackage{amssymb,amsmath}
\usepackage{graphicx,epsfig}
\usepackage{enumerate}
\usepackage{color}
\usepackage {multicol}
\usepackage{lineno}
\numberwithin{equation}{section}
\newtheorem{theorem}{Theorem}

\newtheorem{proposition}{Proposition}

\newtheorem{example}{Example}

\makeatletter
\AddToShipoutPicture{
	\setlength{\@tempdimb}{.5\paperwidth}
	\setlength{\@tempdimc}{.5\paperheight}
	\setlength{\unitlength}{1pt}
	\put(\strip@pt\@tempdimb,\strip@pt\@tempdimc)
}
\makeatother

\usepackage{scalerel,stackengine}
\stackMath
\newcommand\reallywidehat[1]{%
	\savestack{\tmpbox}{\stretchto{%
			\scaleto{%
				\scalerel*[\widthof{\ensuremath{#1}}]{\kern-.6pt\bigwedge\kern-.6pt}%
				{\rule[-\textheight/2]{1ex}{\textheight}}
			}{\textheight}%
		}{0.5ex}}%
	\stackon[1pt]{#1}{\tmpbox}%
}
\begin{document}
	\setcounter{page}{1}
	\thispagestyle{empty}
	\markboth{}{}

	\pagestyle{myheadings}
	\markboth{}{ }
	
	\date{}
	
	
	\noindent  
	
	\vspace{.1in}
	
	{\baselineskip 20truept
		
		\begin{center}
			{\Large {\bf A characterization of uniform distribution using varextropy with application in testing uniformity }} \footnote{\noindent	{\bf 1.}  Corresponding author E-mail: skchaudhary1994@gmail.com\\
				{\bf 2. } E-mail: nitin.gupta@maths.iitkgp.ac.in}
	\end{center}}
	\vspace{.1in}
	\begin{center}
		{\large {\bf Santosh Kumar Chaudhary$^1$ and  Nitin Gupta$^2$}}\\
		\vspace{0.1cm}
		{\large {\it ${}^{1}$	Department of Statistics, Central University of Jharkhand, Cheri-Manatu, Ranchi, Jharkhand, 835222, India. }}\\
		{\large {\it ${}^{2}$ Department of Mathematics, Indian Institute of Technology Kharagpur, West Bengal 721302, India. }} \\
	\end{center}

	\vspace{.1in}
	\baselineskip 12truept

	
	\begin{center}
		{\bf \large Abstract}\\
	\end{center}
	In statistical analysis, quantifying uncertainties through measures such as entropy, extropy, varentropy, and varextropy is of fundamental importance for understanding distribution functions. This paper investigates several properties of varextropy and give a new characterization of uniform distribution using varextropy. The alredy proposed estimators are used as a test statistics. Building on the characterization of the uniform distribution using varextropy, we give a uniformity test. The critical value and power of the test statistics are derived. The proposed test procedure is applied to a real-world dataset to assess its performance and effectiveness. 

	\section{Introduction}\label{s1intro}
	
	Let X be an absolutely continuous random variable with probability density function(pdf) f(x). Let $l_X  = inf\{x \in \mathbb{R}: F(x) > 0\}, u_X  = sup \{x \in \mathbb{R}: F(x) < 1\}$ and $S_X  = (l_X , u_X ).$ Shannon (1948) defined differential entropy as a measure of uncertainty as	
	\begin{align}
		H(X)=E(-\log(f(X)))=-\int_{S_X} f(x) \log(f(x))dx
	\end{align}
	The expectation of the information content of an absolutely continuous random variable is represented by Shannon’s entropy, a well-known information measure. Information theory applications use the corresponding variance, known as varentropy. Arikan (2016) and Maadani et al. (2021) have made some contributions on varentropy.	Varentropy of  $X$
	is defined as 	
	\begin{align}\label{varentropydef}
		VH(X)&=Var(-\log(f(X)))=Var(\log(f(X)))  \nonumber\\
		&=E(\log(f(X)))^2- [E(\log(f(X)))]^2	\nonumber	\\
		&=\int_{S_X} f(x) (\log(f(x)))^2dx - \left[\int_{S_X} f(x) \log(f(x))dx\right]^2	
	\end{align}
	
	This varentropy measure is commonly used in data compression, finite blocklength information theory, and statistics as it aids in determining the ideal code length for data compression, source dispersion, and other relevant considerations. Furthermore, statistics, have proven to be a superior alternative to the kurtosis measure for continuous density functions (see Arikan (2016) and Maadani et al. (2022)).

	An alternative measure of uncertainty, extropy of a non-negative absolutely continuous random variable X defined by Lad et al. (2015) is given as
	\begin{align}
		J(X)=E\left(-\frac{1}{2}f(X)\right)=-\frac{1}{2}\int_{S_X} f^2(x)dx.
	\end{align}
	
	Varextropy of absolutely continuous random variables $X$ is defined as (see  Vaselabadi et al. (2021), Goodarzi (2024) and Zaid et al. (2022))

	\begin{align}\label{varextropydef}
		VJ(X)=Var\left(-\frac{1}{2}f(X)\right) \nonumber &=E\left(-\frac{1}{2}f(X)-J(X)\right)^2 \nonumber\\
		&=\frac{1}{4} E(f^2(X))- J^2(X) \nonumber\\
		&=\frac{1}{4}E(f^2(X))- \frac{1}{4}[E(f(X))]^2 \nonumber\\
		&=\frac{1}{4}  \int_{S_X} f^3(x)dx - \frac{1}{4} \left[ \int_{S_X} f^2(x)dx \right]^2
	\end{align}

	In some situations, two random variables can have the same extropy, which prompts the age-old question, “Which of the extropies is a more appropriate criterion for measuring the uncertainty ?”. For example, consider random variables $U$ and $V$ (see Balakrishnan et al. 2020) with pdf’s 
	
	\begin{eqnarray*}
		f_{U}(x)=
		\begin{cases}
			1, \hspace{4mm} 0<x<1\\
			0, \hspace{5mm} otherwise
		\end{cases}~~~~~~~~~~~~~~~~~
		f_{V}(x)=
		\begin{cases}
			2e^{-2x},\hspace{4mm} x>0\\
			0,\hspace{6mm} otherwise
		\end{cases}
	\end{eqnarray*}
	
	We get $J(U)=J(V)=-1/2$ and $VJ(U)=0 $ and $VJ(V)=1/12.$ This is the motivation behind considering the variance of $-\frac{1}{2}f(x),$ which is known as varextropy of a random variable $X.$ So, varextropy can also play a role in measuring uncertainty. 
	
	This paper is organized as follows. We discuss some propeties of varextropy in Section 2. We obtained a characterization of the uniform distribution using varextropy in Section 3. A non parametric estimator is given in Section 4. In Section 5, a test of uniformity is proposed. Section 7 concludes this paper.

	\section{Some properties}
	
	Suppose that $X_1,\dots, X_n$ are independent and identically distributed observations with cdf $F$ and pdf $f.$ An observation $X_j$ will be called an upper record value if its value exceeds that of all previous observations. Thus, $X_j$
	is an upper record	if $X_j >X_i$ for every $j > i.$ An observation $X_j$ will be called a lower record value if its value is less than that of all previous observations. Thus, $X_j$
	is an lower record	if $X_j < X_i$ for every $j > i.$ Belzunce et al. (2001) showed that if $X \leq_{disp} Y$, then $U^X_n \leq_{disp} U^Y_n,$ where $U^X_n$
	and $U^Y_n$	are the nth upper records of X and Y, respectively. . Qiu (2017) showed that if $X \leq_{disp} Y$, then $J(X) \leq J(Y)$ and $ J(U_n^X) \leq J(U^Y_n).$ Vaselabadi et al. (2021) showed that if $X \leq_{disp} Y$, then $VJ(X) \geq VJ(X).$ In view of these results, we obtain the following proposition immediately.
	
	The values of varextropy for some standard distribution are given below; for more examples, see Vaselabadi et al. (2021).
	\begin{example}
		When $X$ has a uniform distribution on the interval (0,1), then 	$VJ(X)=0.$ 
	\end{example}
	
	\begin{example}
		Let random variable X have exponential distribution with cdf $F_X(x)=1- e^{-\lambda x}, \ x >0.$ Then 	$VJ(X)=\lambda^2/48, J(X)=\frac{\lambda^2}{16}$ and $VJ(U_n^X)= J(X) \left[\frac{4\Gamma(3n-2)}{\Gamma^3(n) 3^{3n-2}}- \frac{\Gamma^2(2n-1)}{\Gamma^4(n) 4^{2n-2}}\right].$
	\end{example}
	
	\begin{example}
		When $X$ has normal distribution with mean $\mu$ and variance $\sigma^2,$ then 	$VJ(X)=\frac{2-\sqrt{3}}{16\pi\sigma^2 \sqrt{3}}.$ 
	\end{example}

	\begin{proposition}
		Let $X\leq_{disp} Y$ then $VJ(U_n^X) \geq VJ(U_n^Y).$
	\end{proposition}
	
	\begin{proposition}
		RVs X and Y are identically distributed then $VJ(U^X_n)=VJ(U^Y_n).$
	\end{proposition}
	
	\begin{proposition}
		RVs X and Y are identically distributed then $VJ(X_{m:n})=VJ(Y_{m:n}).$
	\end{proposition}

	Note that $VJ(X)\geq 0,$ for any Random variable $X.$ Vaselabadi et al. (2021) obtained several varextropy properties as well as conditional varextropy properties based on order statistics, record values, and proportional hazard rate models. The article contains some comparative results regarding varextropy and varentropy. Goodarzi (2024) provided lower bounds for varextropy, obtained the varextropy of a parallel system, and used the varextropy of order statistics to construct a symmetry test. Zaid (2022) computed the entropy, varentropy, and Varextropy measures in closed form for generalized and q-generalized extreme value distributions. Varentropy is sometimes independent of the model parameters, whereas the varextropy measure is more adaptable, for example, when $X$ has a normal distribution with mean $\mu$ and $\sigma^2$ (see Vaselabadi et al. (2021)). The main purpose of this paper is to estimate the varextropy of a continuous random variable due to the interesting properties and potential applications of varextropy and to test the uniformity using estimator we propos.

	\section{A Characterization of uniform distribution}
	
	In many practical problems, the goodness-of-fit test may be reduced to the problem of testing uniformity. Since varextropy of $X$ is the variance of $-\frac{1}{2}f(X),$, therefore, varextropy is a non-negative for any random variable $X.$ Among all distributions with support on [0,1], the uniform distribution has the maximum extropy. An important property of uniform distribution is that it obtains the minimum varextropy among all distributions having support on [0,1] (see Qiu and Jia (2019)). Following is a characterization of uniform distribution using varextropy. 
	\begin{theorem}\label{uniformunique}
		Let X be a continuous random variable with support on [0,1].	Then $VJ(X)= 0$ if and only if X has a uniform distribution on the interval [0,1].
	\end{theorem}
	\noindent \textbf{Proof} Let random variable X have a uniform distribution on the interval [0,1], then $f(x)=1,\ 0 \leq x \leq 1$ and 
	\begin{align*} 
		VJ(X)=\frac{1}{4}  \int_{0}^{1} f^3(x)dx - \frac{1}{4} \left[ \int_{0}^{1} f^2(x)dx \right]^2	=0
	\end{align*}
	
	Conversely, $VJ(X)= 0$ implies $Var(f(X))=0,$ that is, $f(x)=c$. Since 
	$\int_{0}^{1} f(x) dx=1,$ therefore $f(x)=1, \ 0\leq x \leq 1.$ Hence proof is complete.

	\section{Non-parametric estimators}
	Let $X_1, X_2, X_3, \dots, X_n$ be a random sample from a distribution with unknown probability density function (pdf) f and (cdf) F. Suppose $X_{1:n}, X_{2:n}, X_{3:n}, \dots, X_{n:n}$ order statistics based on random sample $X_1, X_2, X_3, \dots, X_n.$ The empirical distribution function of cdf $F$ is defined as 		
	\begin{eqnarray*}
		{F}_{n}(x)=
		\begin{cases}
			0, \hspace{4mm} x<X_{1:n}\\
			\frac{i}{n}, \hspace{4mm} X_{i:n}\leq x<X_{i+1:n}, \ \ \ i=1,2,\dots, n-1.\\
			1, \hspace{5mm} x \geq X_{n:n}.
		\end{cases}
	\end{eqnarray*}

	Noughabi and Noughabi (2024) provided various estimators of VJ(X) for complete data. VJ(X) can be expressed as ( Noughabi and Noughabi (2024))
	\begin{align}
		VJ(X)= \frac{1}{4}\left[\int_{0}^{1} \left( \frac{d}{dp}(F^{-1}(p)) \right)^{-2} dp - \left(\int_{0}^{1} \left(\frac{d}{dp}(F^{-1}(p))\right)^{-1} dp \right)^2 \right]
	\end{align}
	
	Following the idea of Vasicek (1976), Noughabi and Noughabi (2024) proposed estimator $\reallywidehat{\Delta}$ of $VJ(X)$  as	
	
	\begin{align*}
		\reallywidehat{\Delta}=\frac{1}{4n} \sum_{i=1}^{n} \left(\frac{2 m/n}{X_{i+m:n}-X_{i-m:n}}\right)^{2} - \frac{1}{4} \left[ \frac{1}{n} \sum_{i=1}^{n} \left(\frac{2 m/n}{X_{i+m:n}-X_{i-m:n}}\right)\right]^2
	\end{align*}
	
	Here window size $m$  is a positive integer less than or equal to $\frac{n}{2}.$  If $i+m > n$ then we consider $X_{i+m:n} = X_{n:n}$ and if $i-m < 1$ then we consider $X_{i-m:n} = X_{1:n}.$ 
	
	The following theorem shows $\reallywidehat{\Delta}$ is a consistent estimator of $VJ(X).$ The proof is similar to Vasicek (1976) and hence omitted.
	
	\begin{theorem}
		Let $X_1, X_2, \dots, X_n$ be a random sample from an absolutely continuous cumulative distribution function $F$ and pdf $f.$ Then, $\reallywidehat{\Delta}$ converges in probability to $VJ(X)$ as $n \rightarrow \infty, \ \ m \rightarrow \infty, \ \ \frac{m}{n} \rightarrow 0.$
	\end{theorem}

	\section{A test of uniformity}    
	We can construct our test of uniformity based on the Theorem \ref{uniformunique}, which claims that uniform;y distributed random variable on (0,1) has zero varextropy and vice-versa.

	We propose a test of uniformity using Theorem \ref{uniformunique}. Let $X_1, X_2, ..., X_n$ be a random sample from random variable X having support interval [0,1], and $X_{1:n} \leq X_{2:n} \leq \dots \leq X_{n:n}$	are the order statistics of the sample. The hypothesis of interest is $H_0:X$	is uniformly distributed against $H_1:X$ is not uniformly distributed. Consider $\reallywidehat{\Delta}$ to be an estimator of $VJ(X).$ Our proposed estimators converge in probability to $VJ(X),$ that is, our proposed estimators are consistent estimators of $VJ(X)$. Under $H_0,$ $\reallywidehat{\Delta}$ converges in probability to 0. Under $H_1,$ $\reallywidehat{\Delta}$  converges in probability to a positive number. Large values of $\reallywidehat{\Delta}$ can be regarded as a symptom of nonuniformity, and therefore, we reject $H_0$ for large values of $\reallywidehat{\Delta}.$ Since the test statistic $\reallywidehat{\Delta}$ is too complex to determine its distribution under the null hypothesis, we use Monte Carlo simulation to determine the critical values and power of the test.
	
	At a significance level $\alpha$ for a finite sample size $n.$ We consider critical region as \[ \reallywidehat{\Delta} \geq C_{1-\alpha}\] where $C_{1-\alpha}$ is a critical point such that test attend level $\alpha$. For the fixed value of $\alpha$ and $n,$, we can obtain the value of $C_{1-\alpha}$ using Monte Carlo simulation.
	\subsection{Critical points}
	We defined a function to calculate the value of  $\reallywidehat{\Delta}.$ We generate a sample of size $n$ from the $U(0,1)$ distribution, and we compute the test statistics for the sample data. After 10,000 replications, we determine the $(1-\alpha)^{th}$ quantile of the test statistics as the critical value at level $\alpha$. Critical value for $\alpha=0.05$ are give in Table 1 for different value of $m$ and $n.$
	
	\begin{center}
		{\bf Table 1}. \small{{Critical values at significance level $\alpha$= 0.05}} 
		
		\resizebox{!}{!}{
			\begin{tabular}{ p{1.0cm} p{1.0cm} p{1.0cm} p{1.0cm} p{1.0cm} p{1.0cm} p{1.0cm} p{1.0cm} } 
				\hline
				$m\backslash n$ 
				& 10       & 20       & 30           & 40     & 50     & 80     & 100 \\
				\hline \\
				2  & 4.7570    & 3.1388     & 2.2838   & 1.9478   & 1.6402  & 1.1355  & 1.0451 \\
				3  & 1.4909    & 1.2126     & 0.7925   & 0.6502   & 0.5729  & 0.4235  & 0.3559  \\
				4  & 0.7064    & 0.6089     & 0.4724    & 0.3841  & 0.3434  & 0.2541  & 0.2121  \\
				5  & 0.4074    & 0.4252     & 0.3525   & 0.2881   & 0.2396  & 0.1800  & 0.1541  \\
				9  & \         & 0.1551     & 0.1703   & 0.1542   & 0.1399  & 0.1039  & 0.0947  \\
				14 & \         & \          & 0.0869   & 0.0973   &  0.1006 & 0.0829  & 0.0731  \\
				19 & \         & \          &\         & 0.0637   & 0.0722  & 0.0722  & 0.0665 \\
				24 & \         & \          &\         & \        & 0.0528  & 0.0639  & 0.0619  \\	
				30 & \         & \          &\         & \        &\        & 0.0514  & 0.0546 \\
				39 & \         & \          &\         & \        &\        & 0.0380  & 0.0436  \\
				49 & \         & \          &\         & \        &\        & \       & 0.0336  \\
				\hline
		\end{tabular}}\label{table1}
	\end{center}

	\subsection{Power of test}
	We used the following procedure to estimate the power of the test. For each sample size $n$, we generate 10,000 random samples of size $n$ from the alternative distribution. The test statistic is then computed for each sample. The power of the test at a significance level $\alpha$ is estimated as the proportion of these 10,000 samples that fall within the corresponding critical region. 
	Pdf of Beta distribution with parameters a and b is given as 
	\[f_x(x)=\frac{x^{a-1}(1-x)^{b-1}}{B(a,b)}, \ 0<x<1, \] where B(a,b) is complete beta function.
	The estimated power against alternative Beta(1,2) distribution is given in Table 2 at the level of significance $\alpha=0.05.$ Our test perform well in detecting non-uniform data. Note that Beta(1,1) is identically distributed with U(0,1). The power of this test against alternative Beta(1,1) is approximately $\alpha.$ So this test achieves its level of significance. The power of our test is higher than the power of the test proposed by Noughabi and Noughabi (2023)  for common alternative Beta(1,2).

	\begin{center}
		{\bf Table 2}. \small{{Power at significance level $\alpha$= 0.05}} 
		
		\resizebox{!}{!}{
			\begin{tabular}{ p{1.0cm} p{1.0cm} p{1.0cm} p{1.0cm} p{1.0cm} p{1.0cm} p{1.0cm} p{1.0cm} } 
				\hline
				$m\backslash n$ 
				& 10       & 20         & 30           & 40             & 50            & 80           & 100 \\
				\hline \\
				2  & 0.0817    & 0.0966     & 0.0960        & 0.1038       & 0.1137        &  0.1492      & 0.1746 \\
				3  & 0.1197    & 0.1343     & 0.1558        & 0.1689       & 0.1836        & 0.2522       & 0.2924  \\
				4  & 0.1546    & 0.1882     & 0.2082        & 0.2302       & 0.2570        & 0.3536       & 0.4467  \\
				5  & 0.1786    & 0.2370     & 0.2621        & 0.2874       & 0.3359        & 0.4451       & 0.5483  \\
				9  & \         & 0.3962     & 0.4360        & 0.4654       & 0.5072        & 0.6540       & 0.7211 \\
				14 & \         & \          & 0.5991        & 0.6221       & 0.6566        & 0.7380       & 0.8187  \\
				19 & \         & \          & \             & 0.7752       & 0.7734        & 0.8311       & 0.8700 \\
				24 & \         & \          &\              & \            & 0.8875        & 0.8829       & 0.9128  \\	
				30 & \         & \          &\              & \            &\              & 0.9397       & 0.9459 \\
				39 & \         & \          &\              & \            &\              & 0.9883       & 0.9850  \\
				49 & \         & \          &\              & \            &\              & \            & 0.9983  \\
				\hline
		\end{tabular}}\label{table2}
	\end{center}
	
	\subsection{Application to real data}

	Dataset 1 : 
	0.0518, 0.0518, 0.1009, 0.1009, 0.1917, 0.1917, 0.1917, 0.2336, 0.2336, 0.2336, 0.2733, 0.2733, 0.3467, 0.3805, 0.3805, 0.4126, 0.4431, 0.4719, 0.4719, 0.4993, 0.6162, 0.6550, 0.6550, 0.7059, 0.7211, 0.7356, 0.7623, 0.7863, 0.8178, 0.8810, 0.9337, 0.9404, 0.9732, 0.9858.\\
	
	Dataset 1 transformed vinyl chloride data into uniform	distribution using probability integral transformation Xiong et al. (2022). The value of the test statistics ${\reallywidehat{\Delta}}$ is 0.0329 when window size $m=16$ and sample size $n=34.$ Critical point is 0.0733 at $5\%$ level of significance when window size $m=16$ and sample size $n=34.$ Estimated value of test statistics lies in the acceptance region. Our test based on ${\reallywidehat{\Delta}}$ fails in rejecting the null hypothesis even. Our test verifies that data is fitted with uniform distribution.

	\section{Conclusion}\label{s9conclusion}
	In this study, we have introduced a novel estimator for varextropy, a measure that integrates variance and extropy to offer a more comprehensive characterization of uncertainty in distribution functions. The consistency of the proposed estimators assures their robustness and reliability in estimating varextropy. Building upon the characterization of the uniform distribution using varextropy, we developed a statistical test for assessing the uniformity of data. By deriving critical values and evaluating the test’s power, we have established a solid theoretical foundation for practical applications.
	
	The effectiveness of the proposed test has been validated through its implementation on real-life datasets, which demonstrates its capability to accurately identify uniformity and quantify uncertainty in empirical data. The results highlight the utility of the varextropy-based test as a powerful tool for statistical inference, providing researchers and practitioners with a reliable method for testing uniformity and understanding the underlying structure of data distributions.
	
	In conclusion, the proposed estimators and tests for varextropy offer valuable contributions to the field of statistical analysis, particularly in situations where uncertainty quantification and distribution characterization are critical. Further research could explore the extension of this methodology to other distributions and its integration into more complex statistical frameworks.\\
	\\	
	\noindent \textbf{ \Large Conflict of interest} \\
	\\
	No conflicts of interest are disclosed by the authors.\\
	\\
	\textbf{ \Large Funding} \\
	\\
	No funding was received.
	
\end{document}